\documentclass[11pt,a4paper]{amsart}
\usepackage{exscale,amssymb,amsmath,bbm,amsthm}
\usepackage[svgnames]{xcolor}
\usepackage{graphicx}
\usepackage[a4paper,margin=1in,includefoot,footskip=30pt]{geometry}
\usepackage[foot]{amsaddr}
\usepackage{aliascnt}
\definecolor{LinkBlue}{HTML}{0057B8}
\definecolor{CiteBlue}{HTML}{5A3E9B}
\definecolor{UrlRed}{HTML}{B23A2B}
\usepackage[colorlinks=true,linkcolor=LinkBlue,citecolor=CiteBlue,urlcolor=UrlRed]{hyperref}
\usepackage{cleveref}

\newcommand{\height}{\mathrm{height}}
\newcommand{\E}[1]{\ensuremath{\mathbb{E} \left[#1 \right]}}
\newcommand{\Prob}[1]{\ensuremath{\mathbb{P} \left(#1 \right)}}
\newcommand{\var}[1]{\ensuremath{\mathrm{var} \left(#1 \right)}}
\newcommand{\I}[1]{\ensuremath{\mathbbm{1}_{ \{ #1 \} }}}
\newcommand{\R}{\ensuremath{\mathbb{R}}}
\newcommand{\cT}{\mathcal{T}}
\newcommand{\convprob}{\ensuremath{\stackrel{p}{\rightarrow}}}
\newcommand{\Dist}[1]{\mathsf{#1}}
\newcommand{\Exp}{\Dist{Exp}}
\newcommand{\diam}{\ensuremath{\mathrm{diam}}}
\newcommand{\ind}{\ensuremath{\mathbf{1}}}
\newcommand{\eps}{\varepsilon}

\usepackage[sc]{mathpazo}
\usepackage[T1]{fontenc}

\newtheorem{thm}{Theorem}[]
\newaliascnt{lem}{thm}
\newtheorem{lem}[lem]{Lemma}
\aliascntresetthe{lem}
\newaliascnt{defn}{thm}

\aliascntresetthe{defn}
\newaliascnt{prop}{thm}
\newtheorem{prop}[prop]{Proposition}
\aliascntresetthe{prop}
\newaliascnt{claim}{thm}
\newtheorem{claim}[claim]{Claim}
\aliascntresetthe{claim}
\newaliascnt{cor}{thm}
\newtheorem{cor}[cor]{Corollary}
\aliascntresetthe{cor}
\newaliascnt{conj}{thm}

\aliascntresetthe{conj}
\newaliascnt{rem}{thm}

\aliascntresetthe{rem}

\crefname{thm}{Theorem}{Theorems}
\Crefname{thm}{Theorem}{Theorems}
\crefname{lem}{Lemma}{Lemmas}
\Crefname{lem}{Lemma}{Lemmas}
\crefname{defn}{Definition}{Definitions}
\Crefname{defn}{Definition}{Definitions}
\crefname{prop}{Proposition}{Propositions}
\Crefname{prop}{Proposition}{Propositions}
\crefname{claim}{Claim}{Claims}
\Crefname{claim}{Claim}{Claims}
\crefname{cor}{Corollary}{Corollaries}
\Crefname{cor}{Corollary}{Corollaries}
\crefname{conj}{Conjecture}{Conjectures}
\Crefname{conj}{Conjecture}{Conjectures}
\crefname{rem}{Remark}{Remarks}
\Crefname{rem}{Remark}{Remarks}
\crefname{section}{Section}{Sections}
\crefname{subsection}{Section}{Sections}
\crefname{appendix}{Appendix}{Appendices}
\Crefname{appendix}{Appendix}{Appendices}

\title{Cover times of random recursive trees and linear preferential attachment trees}
\author{Omer Angel}
\email{angel@math.ubc.ca}
\thanks{OA was supported in part by NSERC of Canada, SLMath Clay Senior Scholarship, and Magdalen College, Oxford}
\author{Christina Goldschmidt}
\email{christina.goldschmidt@stats.ox.ac.uk}
\author{Peleg Michaeli}
\email{peleg.michaeli@maths.ox.ac.uk}
\thanks{P.M.'s research was supported by UKRI under the Horizon Europe Guarantee for MSCA Postdoctoral Fellowships (EP/Z001781/1, IRIS)}
\date{\today}

\begin{document}

\begin{abstract}
We study the cover times of uniform random recursive trees and linear preferential attachment trees,
including the negative-parameter models corresponding to random $d$-ary increasing trees.
We prove explicit limits in probability and expectation when the cover times are rescaled by $n (\log n)^2$.
In particular, for the uniform random recursive tree the cover time is close to $e n (\log n)^2$ with high probability.
\end{abstract}

\maketitle

A \textbf{recursive tree} on $n$ vertices is a rooted tree labelled by $[n]:= \{1,2,\ldots,n\}$
such that the root has label $1$
and the labels increase along every path away from the root.
We study a sequence of random recursive trees defined as follows.
Fix $a\ge 0$, or $a=-1/d$ for an integer $d\ge 2$.
The tree $T_1$ consists of the single vertex $1$.
For $i \ge 1$, the tree $T_{i+1}$ is obtained from $T_i$ by adding vertex $i+1$ and attaching it to some $j \in [i]$ with probability proportional to $1+a\deg^+_{T_i}(j)$,
where $\deg^+_{T_i}(j)$ is the out-degree (i.e., the number of children) of $j$ in $T_i$.
We write $\cT_a(n)$ for the law of $T_n$.
For $a = 0$, vertex $i+1$ attaches to one of the $i$ earlier vertices, each of which is chosen with probability $1/i$;
this generates the \textbf{uniform random recursive tree} (URRT).
For $a > 0$, this model is known as a \textbf{linear preferential attachment tree}.

This model can be defined also for certain $a<0$.
The construction is well-defined for every $n$ exactly when $a=-1/d$ for an integer $d\ge 1$.
Our results extend naturally to every non-degenerate such case, namely $d\ge 2$.
For $a=-1/d$, a vertex with $k$ children has attachment weight $1-k/d$,
which vanishes when $k=d$.
Equivalently, give each vertex $d$ distinguishable child positions
and place each new vertex in a position chosen uniformly from all those currently vacant.
(This is easily seen to be equivalent to the Eden model on the $d$-ary tree.)
A vertex with $k$ children contributes $d-k$ such positions,
so forgetting the selected positions and keeping only the tree structure gives exactly the attachment probabilities above.
With the positions retained, this gives the random $d$-ary increasing tree studied in \cite[p.~513]{BroutinDevroyeMcLeishDeLaSalle}.
We exclude from this paper the case $a=-1$ in which the tree $T_n$ is not random but is just a path on $n$ vertices.
In this case, the cover time is of order $n^2$ and is not concentrated.
This collection of models has been much studied in the literature;
for background on recursive and preferential-attachment trees, see for example \cite{Pittel,RGCN2Ch8,HolmgrenJanson}.

For a finite connected graph $G$, let $\mathbb{E}_{G,v}$ denote expectation for a discrete-time simple random walk $(X_t)_{t\ge 0}$ on $G$ started at $v$.
For $v \in V(G)$, define the hitting time $\tau_v:=\inf\{t\ge 0:X_t=v\}$, and write $\tau_{\mathrm{cov}}(G) := \max_{v\in V(G)}\tau_v$ for the \textbf{cover time}.
Note that the distribution of the cover time depends, in general, on the starting point of the random walk.
We denote the expectations of the cover and hitting times from the worst starting points by
\[t_{\mathrm{cov}}(G) := \max_{u\in V(G)}\mathbb E_{G,u}[\tau_{\mathrm{cov}}(G)]\] 
and
\[t_{\mathrm{hit}}(G) := \max_{u,v\in V(G)}\mathbb E_{G,u}[\tau_v]. \]

Cover times of finite graphs have been studied extensively since at least the work of Kahn, Linial, Nisan and Saks~\cite{KLNS}.
For random graphs, there is by now a substantial literature,
including results for Erd\H{o}s--R\'enyi graphs~\cite{Jonasson,CooperFriezeSparse,CooperFriezeGiant}, random regular graphs~\cite{CooperFriezeRegular}, random geometric graphs~\cite{CooperFriezeGeom}, and random graphs with given degree sequence~\cite{AbdullahCooperFrieze}.
Specializing to trees, Aldous~\cite{Aldous} studied balanced $d$-ary trees and uniform random labelled trees (Cayley trees);
Brightwell and Winkler~\cite{BrightwellWinkler} showed that, among all choices of an $n$-vertex tree and starting vertex, the minimum expected cover time is attained by a star started from a leaf,
while Feige~\cite{FeigeCoupons} showed that the maximum is attained by a path started from a middle vertex.
More generally, Feige~\cite{FeigeLower} proved that every connected $n$-vertex graph has expected cover time at least $(1+o(1))n\log n$.

The purpose of this note is to find the asymptotic behaviour of $\tau_{\mathrm{cov}}(T_n)$ for the recursive tree models defined above. Let us emphasise that there are two layers of randomness in this model: the underlying randomness of the tree $T_n$ and the randomness associated with the random walk. Some of our results on the cover time will be \emph{quenched}, that is conditional on the underlying tree $T_n$, and others will be \emph{annealed}, that is where both tree and walk are taken to be random. Our primary interest is in the annealed behaviour, but several of the technical tools we use apply only in the quenched setting.
The natural starting point for the random walk in our setting is the root of the tree.

To state our results, let us introduce $W$ for the principal branch of the Lambert $W$ function (sometimes denoted $W_0$). 
This is the unique continuous solution to $x = W(x)e^{W(x)}$ for $x \ge -1/e$ which satisfies $W(0)=0$.
Note that $W$ takes values in $[-1,\infty)$, and is increasing to infinity on its domain (roughly as $\log x$).
We write $o_p(1)$ for any quantity that converges to $0$ in probability as $n\to\infty$.
Our main result can now be stated.

\begin{thm}\label{thm:main}
  Fix $a\ge 0$, or $a=-1/d$ for an integer $d\ge 2$.
  For each $n$, let $T_n\sim\cT_a(n)$,
  and let $\tau_{\mathrm{cov}}(T_n)$ be the cover time of a simple random walk on $T_n$ started at the root. Then, as $n \to \infty$,
  \[
    \tau_{\mathrm{cov}}(T_n)= (1+o_p(1)) \frac{\gamma^2}{(1+a)^2} n (\log n)^2,
  \]
  where 
  \[
    \gamma = \gamma(a) = \sqrt{e} \exp(W(a/(2 \sqrt{e}))) \sqrt{1 + 2 W(a/(2 \sqrt{e}))}.
  \]
  In particular we have $\gamma(0) = \sqrt{e}$ and so for the URRT,
  \[
    \tau_{\mathrm{cov}}(T_n)= ( e + o_p(1)) n (\log n)^2.
  \]
\end{thm}

\begin{figure}
\begin{center}
\includegraphics[width=7cm]{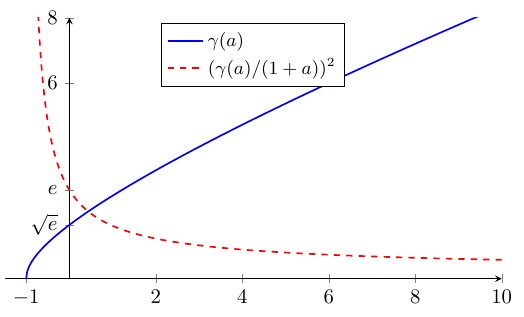}
\end{center}
\caption{The functions $\gamma(a)$ and \mbox{$(\gamma(a)/(1+a))^2$}
from \cref{thm:main}.}
\label{fig:gamma}
\end{figure}

See \Cref{fig:gamma} for an illustration of the relevant constants as functions of $a$.
Note that for $a<0$ not of the form $-d^{-1}$ we do not have any model realizing the specified constants.
Interpolating between $d$-ary graphs, e.g.\ with random degrees would be one option, though we do not pursue this here.

\begin{cor}\label{cor:expectation}
  Under the assumptions of \cref{thm:main}, as $n \to \infty$,
  \[
    \E{t_{\mathrm{cov}}(T_n)} = (1+o(1)) \frac{\gamma^2}{(1+a)^2} n (\log n)^2.
  \]
  In particular for the URRT,
  \[
    \E{t_{\mathrm{cov}}(T_n)} = (e+o(1)) n (\log n)^2.
  \]
  Here the expectation is over $T_n$.
\end{cor}

Before proceeding, let us briefly outline the idea behind \cref{thm:main}.
We start with a heuristic analysis.
For clarity, we focus on the URRT at first.
The proof is based on the well-known connection between cover times and Gaussian free fields, see \cite{DingLeePeres,Ding} and \cref{thm:Ding} below.
To apply this we need to estimate the maximum of a GFF on $T_n$.
Now, by the Sudakov--Fernique \cite{Sudakov71,Fernique75} inequality, the expected maximum is at most the expected maximum of independent Gaussian variables
indexed by the vertices, with each variable's variance equal to the variance of the GFF value at that vertex.
The number of vertices in $T_n$ at level $s\log n$ is known~\cite{FuchsHwangNeininger06} to be $n^{s-s\log(s)+o(1)}$, and the maximum of that many i.i.d.\ $N(0,s\log n)$ variables is concentrated around $\sqrt{s\log n}\cdot\sqrt{2(s-s\log s)\log n} = s\sqrt{2(1-\log s)} \log n$.
This is maximized at $s=\sqrt{e}$, and so the maximum of the GFF is at most $(\sqrt{e} + o_p(1))\log n$, heuristically attained at some vertices of height $\sqrt{e}\log n$.
It is possible to get a matching lower bound (which we do).
Typically two vertices of $T_n$ have a nearest common ancestor which is quite close to the root.
This implies that the values of the GFF are not too far from independent.

Our actual proof follows a different path:
we construct the tree $T_n$ simultaneously with the GFF on it, and observe that the maximum of the GFF can be viewed as the maximal displacement of a certain branching random walk.
This in turn can be studied through the general associated theory developed by Biggins \cite{Biggins95,Biggins97}.

\section{Background}
\subsection{Properties of linear preferential attachment trees}

For a tree $T$,
we write $\height(T)$ for the maximum graph distance from the root to a vertex of $T$ and $\diam(T)$ for the diameter of $T$ (that is, the largest graph distance between any two of its vertices).
Then
\[
  \height(T)\le \diam(T)\le 2\cdot\height(T).
\]
For preferential attachment trees, 
the first-order height asymptotic was proved for $a\ge0$ by Pittel~\cite[Theorem~1]{Pittel}, 
and for $a=-1/d$ by Broutin et al.~\cite[Theorem~3]{BroutinDevroyeMcLeishDeLaSalle}.
We use below the following result giving convergence in $L^1$.

\begin{thm}[\cite{HolmgrenJanson}, Theorem~13.19 and Example~13.9]\label{thm:height}
  Suppose that $a\ge 0$, or that $a=-1/d$ for an integer $d\ge 2$.
  Then
  \[
    \frac{\height(T_n)}{\log n}\longrightarrow\frac{1}{(1+a)\xi}
    \qquad\text{in }L^1,
  \]
  where $\xi$ is the unique root in $(0,1)$ of the equation $a \xi + \log \xi + 1 = 0$.
\end{thm}

\subsection{Continuous time construction}

A well-known tool for studying random recursive trees and preferential attachment trees is their continuous-time construction as Crump--Mode--Jagers branching processes;
see Pittel~\cite{Pittel} for $a\ge 0$, and Holmgren and Janson~\cite{HolmgrenJanson} for the $d$-ary case.
We recap this construction for completeness.
Throughout, for $\lambda\ge 0$, let $\Exp(\lambda)$ denote the exponential distribution with rate $\lambda$.
(Here we take $\Exp(0)=\infty$ with probability 1.)
We embed the sequence $(T_n)_{n\geq 1}$ in continuous time as follows.
Consider a continuous-time pure birth process that generates a sequence of growing trees.
For $t\ge 0$, denote by $T(t)$ the tree at time $t$.
The process starts at $t=0$ with $T(0)$ consisting of a single root vertex,
labelled $1$.
Each vertex $j$ in the tree is endowed with an independent \textbf{exponential clock}.
At time $t$, the rate of its clock is $1+a\cdot \deg_{T(t)}^+(j)$.
When the clock belonging to vertex $j$ rings, a new vertex (with the next free label) is created,
and attached to $j$ as a child.
The rate of $j$'s clock is then updated, and an independent clock is started for the new vertex.
Let $\beta_k$ be the random time at which vertex $k$ first appears in the tree.

\begin{claim}\label{cl:beta-sum}
  Set $E_i:=\beta_{i+1}-\beta_i$ for $i\ge 1$.
  Then the random variables $(E_i)_{i\ge 1}$ are independent, with $E_i\sim\Exp((1+a)i-a)$.
  Moreover, the embedded chain $(T(\beta_n))_{n\ge 1}$ has the same joint law as $(T_n)_{n\ge 1}$.
\end{claim}

\begin{proof}
  Fix $i\ge 1$ and consider the process at time $\beta_i$,
  when $T(\beta_i)$ has vertex set $[i]$.
  Conditional on the evolution up to this time,
  the times from $\beta_i$ until the clocks next ring are, by the memoryless property, independent exponential random variables with respective rates $1+a\deg^+_{T(\beta_i)}(v)$, for $v\in[i]$.
  Since the out-degrees in $T(\beta_i)$ sum to $i-1$, the total clock rate is $(1+a)i-a$.
  Hence their minimum $E_i=\beta_{i+1}-\beta_i$ has distribution $\Exp((1+a)i-a)$,
  while vertex $i+1$ attaches to $v\in[i]$ with probability proportional to $1+a\deg^+_{T(\beta_i)}(v)$.
  Since the total clock rate depends only on $i$,
  the conditional law of $E_i$ does not depend on the evolution up to time $\beta_i$.
  In particular, $E_i$ is independent of $E_1,\ldots,E_{i-1}$,
  so the variables $(E_i)_{i\ge 1}$ are independent.
  Since $T(\beta_1)=T_1$,
  induction shows that the embedded chain has the same joint law as $(T_n)_{n\ge 1}$.
  Finally, $\beta_1=0$, so $\beta_n=\sum_{i=1}^{n-1}E_i$.
\end{proof}

\subsection{Gaussian free fields on trees} \label{subsec:GFF}

The Gaussian free field (GFF) on graphs has well known connections to random walks and cover times.
See e.g.\ \cite{LyonsPeres} for an introduction to discrete GFFs.
The GFF on a tree $T$, pinned at a vertex $v_0$, has a particularly simple description.
Root $T$ at $v_0$ and set $\eta_{v_0}=0$.
Then, recursively, whenever $v$ is the parent of $u$ in the tree,
let $\eta_u-\eta_v$ have a standard normal distribution,
with these increments independent over the edges of $T$.

For general finite graphs, Ding, Lee and Peres~\cite{DingLeePeres} proved that the maximal expected cover time over starting vertices is, up to universal constant factors,
the number of edges times the square of the expected maximum of the GFF.
We use the following sharper tree-specific concentration estimates of Ding~\cite{Ding}.
Ding states them for continuous-time random walks, but notes after his Corollary~1.3 that they also hold in discrete time.

\begin{thm}[\cite{Ding}, Theorem~1.2]\label{thm:Ding}
  There exist universal constants $C,c>0$ such that the following holds.
  Let $T=(V,E)$ be a tree with at least two vertices and let $v_0\in V$.
  Let $\{\eta_v\}_{v\in V}$ be the GFF on $T$ pinned at $v_0$.
  Then, for the random walk started at $v_0$ and every $\lambda\ge 1$,
  \[
    \Prob{\left|\tau_{\mathrm{cov}}(T)
                - |E| \left( \E{\sup_{v \in V} \eta_v}\right)^2 \right|
          \ge \lambda |E| \sqrt{\diam(T)}\E{\sup_{v \in V} \eta_v}}
    \le C e^{-c \lambda}.
  \]
\end{thm}

\begin{cor}[\cite{Ding}, Corollary~1.3]\label{cor:Ding}
  There exist universal constants $C,c>0$ such that, for every finite rooted tree $T$ with at least two vertices and every $\lambda\ge 0$,
  \[
    \Prob{\left|\tau_{\mathrm{cov}}(T)-t_{\mathrm{cov}}(T)\right|
      \ge \lambda\sqrt{t_{\mathrm{cov}}(T)t_{\mathrm{hit}}(T)}}
    \le Ce^{-c\lambda}.
  \]
  Here the random walk is started at the root.
\end{cor}

\section{Proofs}
\subsection{Branching random walks}

We augment the continuous-time tree-growth process by assigning a spatial position to each vertex, thereby constructing the GFF on the tree at the same time.
This defines a \textbf{branching random walk (BRW)}.
The root (vertex $1$) is located at position $\eta_1=0$.
Then, recursively, when a new vertex $j$ is born to a parent $i$ at time $\beta_j$,
its position is set to be $\eta_j=\eta_i+N_j$,
where the $\{N_j\}_{j\ge 2}$ are independent and identically distributed standard normal random variables.
The sequence $(N_j)_{j\ge 2}$ is independent of the tree-growth process.
The joint process --- comprising the tree's evolution and the assignment of spatial positions --- is a specific instance of a general continuous-time branching random walk model which has been studied extensively, including in seminal works of Biggins~\cite{Biggins95,Biggins97}.
Evidently, the set of positions $\{\eta_i\}_{i\in[n]}$ has the same law as the Gaussian free field on $T_n$ pinned at vertex $1$.
Therefore, the maximum of the GFF corresponds to the rightmost position of a particle in this BRW.

Note that if instead of Gaussian displacements we take a deterministic $+1$ displacement between a particle and its parent, the position of a particle is the height of the corresponding vertex in the tree.
This has been used in the past to study the height of the tree (see, e.g.,~\cite{HolmgrenJanson}).

In order to relate our statements to the notation and terminology used by Biggins and Pittel, we let
the maximal possible number of children of a vertex be denoted by
\[
  K_a:=
  \begin{cases}
    \infty,&a\ge 0,\\
    d,&a=-1/d.
  \end{cases}
\]
The offspring point process for a vertex gives the spatial displacement and relative birth times of its children.
In our case, this is
\[
  Z=\{(N_k,S_k)\}_{k=1}^{K_a},
\]
where $\{N_k\}_{1\le k\le K_a}$ are independent standard normal random variables.
The birth times $S_k$ are defined in terms of 
$\{X_k\}_{1\le k\le K_a}$ which are independent exponential random variables with $X_k\sim\Exp(1+(k-1)a)$, and $S_k=\sum_{i=1}^k X_i$.
For $t\ge 0$ and $x\in\R$,
let $P_t(x)$ be the number of particles present at time $t$ whose positions are at least $x$. The rightmost particle is given by
\[
  R_t=\sup\{x\in\R\ :\ P_t(x)\ge 1\}.
\]

The analysis of the branching random walk relies on the two-dimensional Laplace transform of the intensity measure of $Z$
on the domain $\R\times(a,\infty)$, which is defined by
\[
  m(\theta, \phi) := \E{\sum_{k=1}^{K_a} e^{-\theta N_k - \phi S_k}}.
\]

Following Biggins, define for this BRW,
\[
  \alpha(\theta)=\inf\{\phi\ :\ m(\theta,\phi)\le 1\}
  \qquad\text{and}\qquad
  \alpha^*(x)=\inf_{\theta<0}\{x\theta+\alpha(\theta)\}.
\]
Set $\gamma_0=\inf\{x:\alpha^*(x)<0\}$.
The following consequence of Biggins~\cite[Corollary 2.2]{Biggins97} shows that $R_t/t$ converges almost surely to $\gamma_0$.
Biggins' corollary is stated for a more general model;
in \cref{app:biggins} we show how our process fits that model and verify the required hypotheses.

\begin{prop}[Right-edge theorem]\label{prop:biggins-right}
  For the branching random walk constructed above, we have
  \[
    \frac{R_t}{t}\to\gamma_0
  \]
  almost surely as $t\to\infty$.
\end{prop}

Before proceeding, let us give some intuition for this result, which may seem somewhat obscure to readers not familiar with large deviations for branching random walks.
For a vertex $v$, let $p(v)$ denote its position.
Write $A_t$ for the set of vertices born by time t, and define $w_t(\theta) = \E{\sum_{v \in A_t} e^{-\theta p(v)}}$.
It is generally the case for branching random walks that such quantities grow (or decay) exponentially with $t$.
Let us posit that this occurs at some rate $\alpha(\theta)$, i.e.
\[
w_t(\theta) \approx e^{\alpha(\theta)t}.
\]
Note that the root is at position 0.
Condition on the times of birth and relative positions of the offspring of the root individual, and note that the $k$th child of the root, which is born at time $S_k$, is the root of a new, equally distributed branching process.
The time and spatial shift from the root imply this sub-process contributes $w_{t-S_k}(\theta)$ on average, shifted by $e^{-\theta N_k}$ to account for its starting position.
We obtain the equation
\[
w_t(\theta) = 1 + \E{ \sum_{k \ge 1} w_{t-S_k}(\theta) e^{-\theta N_k} \I{S_k \le t}}.
\]
Using the approximate exponential growth, we then get
\[
e^{\alpha(\theta) t} \approx \E{ \sum_{k \ge 1} e^{\alpha(\theta)(t - S_k)} e^{-\theta N_k}}
\]
for large $t$. Cancelling a factor of $e^{\alpha(\theta)t}$ from both sides yields that $\alpha(\theta)$ should be chosen so that
\[
1 = \E{\sum_{k\ge 1} e^{- \alpha(\theta) S_k - \theta N_k}} = m(\theta, \alpha(\theta)).
\]
This gives rise to the definition of $\alpha(\theta)$ given above.

Now consider $P_t(tx) = |\{v \in A_t: p(v) \ge tx\}| \approx  |\{v \in A_t: p(v) \approx  tx\}|$.
Then, for $\theta<0$,
\[
\E{P_t(tx) e^{-\theta t x}} \lesssim w_t(\theta) \approx e^{\alpha(\theta) t}
\]
and so
\[
\E{P_t(tx)} \lesssim \exp \left((\alpha(\theta) + \theta  x) t \right).
\]
Since this holds for all $\theta<0$,
we may optimise to get a tighter bound on the right-hand side by taking $\alpha^*(x) = \inf_{\theta < 0} (\alpha(\theta) + \theta x)$.
A large deviations argument shows that this, in fact, captures the correct exponential rate, so that indeed
\[
\E{P_t(tx)} \approx e^{\alpha^*(x) t}.
\]
The right edge of the population at large time $t$ is then to be found approximately at the point $tx$ where $x$ is such that $\alpha^*(x) = 0$ (i.e.\ the point where $\E{P_t(tx)}$ tips from growing exponentially to decaying exponentially). This is precisely $\gamma_0$.

We now compute these quantities explicitly in our setting:
$m$ in \cref{cl:m}, $\alpha$ and $\alpha^*$ in \cref{cor:alpha},
and finally $\gamma_0$ in \cref{cl:gamma0}.

\begin{claim}\label{cl:m}
  For $\theta\in\R$ and $\phi>a$,
  \[
    m(\theta,\phi)=\frac{1}{\phi-a}\cdot e^{\theta^2/2}.
  \]
\end{claim}

\begin{proof}
  By linearity of expectation and by independence,
  \[
    m(\theta, \phi) = \sum_{k=1}^{K_a} \left(\E{ e^{-\theta N_k}} \cdot \E{e^{-\phi S_k} }\right).
  \]
  The law of $N_k$ does not depend on $k$, and $\E{e^{-\theta N_k}} = \exp(\theta^2/2)$.
  The remaining identity
  \[
    \sum_{k=1}^{K_a}\E{e^{-\phi S_k}}=\frac{1}{\phi-a}
  \]
  is implicit, when $a\ge 0$, in \cite[proof of Theorem~1, p.~340]{Pittel};
  for completeness, and to include the case $a=-1/d$, we derive it explicitly here.
  If $a=0$, then $X_i\sim\Exp(1)$ for every $i$, so
  \[
    Y_k:=\E{e^{-\phi S_k}}
    = \prod_{i=1}^k \E{e^{-\phi X_i}}
    = (1+\phi)^{-k}.
  \]
  Therefore
  \[
    \sum_{k=1}^{\infty}Y_k
    = \sum_{k=1}^{\infty}(1+\phi)^{-k}
    = \frac{1}{\phi},
  \]
  which is the required formula when $a=0$.

  Assume henceforth that $a\ne 0$.
  Letting $(x)^{(k)} = \prod_{i=0}^{k-1} (x+i)$
  and defining $A=1/a$ and $B=(1+\phi)/a$, we have
  \[
    Y_k := \E{e^{-\phi S_k}}
    = \prod_{i=1}^k \E{e^{-\phi X_i}}
    = \prod_{i=1}^k \frac{1+(i-1)a}{1+(i-1)a+\phi}
    = \prod_{i=0}^{k-1} \frac{A+i}{B+i}
    = \frac{A^{(k)}}{B^{(k)}}.
  \]
  If $a>0$, then $K_a=\infty$ and $B-A-1=(\phi-a)/a>0$,
  so Gauss' hypergeometric summation formula gives
  \[
    \sum_{k=0}^{\infty}\frac{A^{(k)}}{B^{(k)}}
    =\frac{B-1}{B-A-1}.
  \]
  If $a=-1/d$, then $A=-d$ and $K_a=d$.
  Moreover, $B^{(k)}\ne 0$ for $0\le k\le d$,
  since $\phi>-1/d$.
  The terminating form of Gauss' formula therefore gives
  \[
    \sum_{k=0}^{d}\frac{A^{(k)}}{B^{(k)}}
    =\frac{B-1}{B-A-1}.
  \]
  Thus, in either case,
  \[
    \sum_{k=0}^{K_a}\frac{A^{(k)}}{B^{(k)}}
    =\frac{B-1}{B-A-1}.
  \]
  Subtracting the $k=0$ term, which is equal to $1$, we obtain
  \[
    \sum_{k=1}^{K_a}Y_k
    =\frac{B-1}{B-A-1}-1
    =\frac{A}{B-A-1}
    =\frac{1}{\phi-a},
  \]
  and the statement follows.
\end{proof}

\begin{cor}\label{cor:alpha}
  For all $\theta\in\R$,
  $\alpha(\theta)=a+e^{\theta^2/2}$,
  and for all $x>0$,
  \[
    \alpha^*(x)
    =a+\exp\left(\frac{1}{2} W(x^2)\right)-x\sqrt{W(x^2)}
    =a+\left(1-W(x^2)\right)\exp\left(\frac{1}{2} W(x^2)\right),
  \]
  where $W$ is the principal branch of the Lambert $W$ function.
\end{cor}

\begin{proof}
  By \cref{cl:m}, for fixed $\theta$ and $\phi>a$ we have $m(\theta,\phi)=e^{\theta^2/2}/(\phi-a)$,
  which is strictly decreasing in $\phi$ with $m(\theta,\phi)\to 0$ as $\phi\to\infty$ and $m(\theta,\phi)\to\infty$ as $\phi\to a$.
  Thus $m(\theta,\phi)\le 1$ iff $\phi-a\ge e^{\theta^2/2}$,
  and the infimum is attained at equality,
  giving $\alpha(\theta)=a+e^{\theta^2/2}$.

  Fix $x>0$ and consider $f(\theta)=x\theta+\alpha(\theta)=x\theta+a+e^{\theta^2/2}$ for $\theta<0$.
  The function $f$ is strictly convex since $f''(\theta)=(1+\theta^2)e^{\theta^2/2}>0$,
  so the minimum over $(-\infty,0)$ is attained at the unique point $\theta_*$ where $f'(\theta_*)=0$.
  We find that point:
  \[
    0=f'(\theta_*)=x+\theta_* e^{\theta_*^2/2}
    \quad\Longleftrightarrow\quad
    \theta_* e^{\theta_*^2/2}=-x.
  \]
  Let $y=\theta_*^2>0$.
  Squaring gives $x^2=ye^{y}$, hence $y=W(x^2)$.
  From $x=-\theta_* e^{y/2}$ and $\theta_*<0$ we get $e^{y/2}=x/\sqrt{y}$.
  Evaluating $f$ at the minimiser,
  \[
    \alpha^*(x)=
    a+e^{y/2}-x\sqrt{y}
    = a+e^{y/2}-ye^{y/2}
    = a+(1-y)e^{y/2}.\qedhere
  \]
\end{proof}

\begin{claim}\label{cl:gamma0}
  Let $w=W(a/(2\sqrt{e}))$.
  Then
  \[
    \gamma_0
    = \sqrt{e}\cdot\exp(w)\cdot\sqrt{1+2w}.
  \]
  In particular, $\gamma_0=\gamma(a)$.
\end{claim}

\begin{proof}
  Writing $y=W(x^2)$, we know by \cref{cor:alpha} that $\alpha^*(x)=a+(1-y)e^{y/2}$.
  Thus we have $\alpha^*(x)<0$ if and only if $(y-1)e^{y/2}>a$.
  Now, $\alpha^*(x)$ is strictly decreasing in $y$ on $(-1,\infty)$.
  Thus there is a unique $y_*$ solving the equation $(y-1)e^{y/2}=a$.
  By the previous equivalence, $\alpha^*(\gamma_0)=0$ if and only if $y_*=W(\gamma_0^2)$.
  Writing $y_*=2z+1$, we have $ze^z=a/(2\sqrt{e})$,
  so $z=w$, and $y_*=1+2w$.
  This implies 
  $W(\gamma_0^2)=1+2w$,
  hence $\gamma_0^2 = (1+2w)e^{1+2w}$,
  from which the statement follows.
\end{proof}

By \cref{prop:biggins-right} and \cref{cl:gamma0}, we obtain the following.

\begin{prop}\label{prop:Rt}
  The position of the rightmost particle satisfies
  $R_t/t\to\gamma(a)$ almost surely, with $\gamma(a)$ as in \cref{thm:main}.
\end{prop}

\begin{proof}
  By \cref{prop:biggins-right},
  \[
    \frac{R_t}{t}\longrightarrow \gamma_0
  \]
  almost surely.
  Finally, by \cref{cl:gamma0} we have $\gamma_0=\gamma(a)$, so $R_t/t\to\gamma(a)$ almost surely.
\end{proof}

We will need the following law of large numbers for $\beta_n$ (the time at which the growing tree first has $n$ vertices).
This lets us apply \cref{prop:Rt} at the random time $t=\beta_n$.

\begin{lem}\label{cl:beta-lln}
  As $n\to\infty$,
  \[
    \frac{\beta_n}{\log n}\to \frac{1}{1+a}
  \]
  almost surely.
\end{lem}

We follow the proof of \cite[Theorem 1]{Pittel}:
First and second moments give concentration of $\beta_n$ near $\frac{\log n}{1+a}$; monotonicity and restriction to a lacunary sequence strengthen to a.s. convergence.

\begin{proof}
  By \cref{cl:beta-sum},
  $\beta_n=\sum_{k=1}^{n-1}E_k$ with independent $E_k\sim\Exp((1+a)k-a)$.
  Thus,
  \[
    \E{\beta_n}
    =\sum_{k=1}^{n-1}\E{E_k}
    =\frac{1}{1+a}\sum_{k=1}^{n-1}\frac{1}{k-\frac{a}{1+a}}
    =\frac{1}{1+a}\log n+O(1),
  \]
  and
  \[
    \var{\beta_n}
    =\sum_{k=1}^{n-1}\var{E_k}
    =\sum_{k=1}^{n-1}\frac{1}{((1+a)k-a)^2}=O(1).
  \]
  Let $n_m=2^m$. For any fixed $\varepsilon>0$, Chebyshev gives
  \[
    \Prob{\left|\beta_{n_m}-\E{\beta_{n_m}}\right|>\varepsilon\log n_m}
    \le \frac{\var{\beta_{n_m}}}{\varepsilon^2\log^2 n_m}
    \le \frac{C}{\varepsilon^2 m^2},
  \]
  for some constant $C<\infty$.
  Since $\sum_{m\ge 1} C/(\varepsilon^2 m^2)<\infty$, the Borel--Cantelli lemma implies
    $\beta_{2^m}/(m\log 2)\to 1/(1+a)$ almost surely.
  For $2^m\le n<2^{m+1}$ we have $\beta_{2^m}\le \beta_n\le \beta_{2^{m+1}}$ and $m\log 2\le \log n\le (m+1)\log 2$, hence
  \[
    \frac{\beta_{2^m}}{(m+1)\log 2}\ \le\ \frac{\beta_n}{\log n}\ \le\ \frac{\beta_{2^{m+1}}}{m\log 2},
  \]
  and both bounds tend to $1/(1+a)$ almost surely.
\end{proof}

The next statement follows easily.
\begin{cor}\label{cor:max-as}
  The GFF on the growing tree satisfies 
  \[
  \frac{\sup_{v\in V(T_n)}\eta_v}{\log n}\to \frac{\gamma}{1+a}
  \]
  almost surely as $n\to\infty$.
\end{cor}

\begin{proof}
  By \cref{prop:Rt}, $R_t/t\to\gamma$ almost surely.
  By \cref{cl:beta-lln}, $\beta_n/\log n\to 1/(1+a)$ almost surely.
  On the intersection of these probability-1 events,
  \[
    \frac{\sup_{v\in V(T_n)}\eta_v}{\log n}
    =\frac{R_{\beta_n}}{\log n}
    =\Big(\frac{R_{\beta_n}}{\beta_n}\Big)\Big(\frac{\beta_n}{\log n}\Big)
    \to \gamma\cdot\frac{1}{1+a}.\qedhere
  \]
\end{proof}

\subsection{From annealed to quenched}

\Cref{cor:max-as} establishes that the maximum of the GFF converges on rescaling when both the tree and the GFF are treated as random.
In order to apply Ding's results from \cref{subsec:GFF}, however, we need the analogous result conditioned on the tree.

The Borell inequality below (also discovered by Tsirelson, Ibragimov, and Sudakov) appears, for example, in \cite[Lemma~4]{Talagrand88}.
\begin{lem}[Borell inequality]\label{lem:talagrand}
  Let $(X_i)_{i\in I}$ be a finite centred Gaussian process, and set
  \[
    M:=\sup_{i\in I} X_i
    \qquad\text{and}\qquad
    \theta:=\sup_{i\in I}\sqrt{\var{X_i}}.
  \]
  Suppose that $\theta>0$, and let $w\in\R$ satisfy
  $\Prob{M\ge w}\le 1/2$. Then, for every $s\ge 0$,
  \[
    \Prob{M\ge w+s}\le \Prob{Y\ge s/\theta},
  \]
  where $Y$ is a standard normal random variable.
\end{lem}

\begin{lem}\label{lem:borell-tree}
  Let $T$ be a rooted tree with at least three vertices and height $h$,
  let $(\eta_v)_{v\in V(T)}$ be the pinned GFF on $T$,
  and set
  \[
    M(T):=\sup_{v\in V(T)}\eta_v.
  \]
  Let $m_T$ be a median of $M(T)$. Then for every $s\ge 0$,
  \[
    \Prob{M(T)\ge m_T+s}\le \Prob{Y\ge s/\sqrt h},
  \]
  where $Y$ is a standard normal random variable.
  In particular, there is a universal constant $C<\infty$ such that
  \[
    \E{(M(T)-m_T)^+}\le C\sqrt h.
  \]
\end{lem}

\begin{proof}
  Since $T$ has at least three vertices, $\Prob{M(T)=0}<1/2$
  (consider any two distinct non-root vertices), while $M(T)$ has no atoms on $(0,\infty)$.
  Hence the median $m_T$ is well defined, positive and $\Prob{M(T)\ge m_T}=1/2$.

  The process $(\eta_v)_{v\in V(T)}$ is a discrete centred Gaussian process.
  For each $v$,
  \[
    \var{\eta_v}=d(\rho,v)\le h,
  \]
  where $d(\rho,v)$ is the graph distance from the root.
  Hence,
  \[
    \sup_{v\in V(T)} \sqrt{\var{\eta_v}}\le \sqrt h.
  \]
  By \cref{lem:talagrand}, applied to the Gaussian process $(\eta_v)_{v\in V(T)}$,
  \[
    \Prob{M(T)\ge m_T+s}\le \Prob{Y\ge s/\sqrt h}
  \]
  for every $s\ge 0$.
  Integrating this bound yields
  \[
    \E{(M(T)-m_T)^+}
    = \int_0^\infty \Prob{M(T)\ge m_T+s}\,ds
    \le \sqrt h \int_0^\infty \Prob{Y\ge u}\,du
    \le C\sqrt h
  \]
  for a universal constant $C$.
\end{proof}

For brevity, set
\[
  \kappa=\kappa(a):=\frac{\gamma(a)}{1+a}.
\]
The next lemma is a quenched concentration result for the maximum of the GFF on $T_n$.

\begin{lem}\label{lem:good-trees}
  There exists $C=C(a)<\infty$ such that, for every fixed $\eps\in(0,1/2)$,
  with probability tending to $1$ as $n\to\infty$,
  the tree $T_n$ satisfies
  \[
    \height(T_n)\le C\log n,
  \]
  and the maximum of the GFF satisfies
  \[
    \Prob{\left|M(T_n)-\kappa\log n\right|>\eps\log n\ \Big|\ T_n}\le \eps.
  \]
\end{lem}

\begin{proof}
  By \cref{thm:height}, there exists $C=C(a)$ such that
  \[
    \Prob{\height(T_n)\le C\log n}\to 1.
  \]
  Let
  \[
    q_n(T):=\Prob{\left|M(T)-\kappa\log n\right|>\eps\log n}.
  \]
  Under the joint law of $(T_n,\eta)$, \cref{cor:max-as} implies
  \[
    \Prob{\left|M(T_n)-\kappa\log n\right|>\eps\log n}\to 0.
  \]
  By definition,
  \[
    q_n(T_n)
    =
    \Prob{\left|M(T_n)-\kappa\log n\right|>\eps\log n\ \Big|\ T_n}.
  \]
  Taking expectations and using the tower property gives
  \[
    \E{q_n(T_n)}
    =
    \Prob{\left|M(T_n)-\kappa\log n\right|>\eps\log n}.
  \]
  So $q_n(T_n)\to 0$ in probability.
  In particular,
  \[
    \Prob{q_n(T_n)\le \eps}\to 1.
  \]
  Combining this with the height bound completes the proof.
\end{proof}

\begin{lem}\label{lem:expectation}
  With conditional expectation given $T_n$,
  \[
    \mathbb E_{T_n}[M(T_n)]
    =
    \left(\kappa+o_p(1)\right)\log n.
  \]
\end{lem}

\begin{proof}
  Write
  \[
    M_n:=M(T_n).
  \]
  Let $m_n$ be any median of the conditional law of $M_n$ given $T_n$.
  Fix $\eps\in(0,1/2)$, and let $C=C(a)$ be as in \cref{lem:good-trees}.
  Let $G_{n,\eps}$ be the event that
  \[
    \height(T_n)\le C\log n
  \]
  and
  \[
    \Prob{\left|M_n-\kappa\log n\right|>\eps\log n\ \Big|\ T_n}\le \eps.
  \]
  On $G_{n,\eps}$, the second condition implies
  \[
    \Prob{M_n\le (\kappa-\eps)\log n\mid T_n}\le \eps
    \qquad\text{and}\qquad
    \Prob{M_n\ge (\kappa+\eps)\log n\mid T_n}\le \eps.
  \]
  Since $\eps<1/2$, this implies
  \begin{equation}\label{eq:mn-close}
    |m_n-\kappa\log n|\le \eps\log n.
  \end{equation}
  By \cref{lem:borell-tree}, for some universal constant $C'>0$,
  \begin{equation}\label{eq:positive-part}
    \mathbb E_{T_n}[(M_n-m_n)^+]\le C'\sqrt{\height(T_n)}
    \le C'\sqrt{C\log n}
  \end{equation}
  on $G_{n,\eps}$.
  Let
  \[
    r_n:=2\eps\log n.
  \]
  Since $M_n\ge 0$, we also have
  \[
    (m_n-M_n)^+\le r_n + m_n\,\I{M_n\le m_n-r_n}.
  \]
  Therefore, on $G_{n,\eps}$,
  \[
    \mathbb E_{T_n}[(m_n-M_n)^+]
    \le 2\eps\log n + m_n\,\Prob{M_n\le m_n-2\eps\log n\mid T_n}.
  \]
  Since $m_n\le (\kappa+\eps)\log n$ and
  \[
    m_n-2\eps\log n\le (\kappa-\eps)\log n,
  \]
  we obtain on $G_{n,\eps}$ that
  \begin{equation}\label{eq:negative-part}
    \mathbb E_{T_n}[(m_n-M_n)^+]
    \le 2\eps\log n + (\kappa+\eps)\eps\log n
    = (\kappa+\eps + 2)\eps\log n.
  \end{equation}
  By \eqref{eq:mn-close}, \eqref{eq:positive-part}, and \eqref{eq:negative-part},
  we find that on $G_{n,\eps}$,
  \begin{equation}\label{eq:mean-close}
    \begin{aligned}
      \left|\mathbb E_{T_n}[M_n]-\kappa\log n\right|
      &\le \left|\mathbb E_{T_n}[M_n]-m_n\right| + \left|m_n-\kappa\log n\right| \\
      &\le \mathbb E_{T_n}[(M_n-m_n)^+] + \mathbb E_{T_n}[(m_n-M_n)^+] + \eps\log n \\
      &\le C'\sqrt{C\log n} + (\kappa+\eps+2)\eps\log n + \eps\log n\\
      &= C'\sqrt{C\log n} + (\kappa+\eps+3)\eps\log n.
    \end{aligned}
  \end{equation}
  Given $\delta>0$, choose $\eps\in(0,1/2)$ so small that
  $
    (\kappa+\eps+3)\eps<\delta/2
  $.
  By \cref{lem:good-trees}, $\Prob{G_{n,\eps}}\to 1$.
  Since $\sqrt{\log n}=o(\log n)$,
  \eqref{eq:mean-close} implies
  \[
    \Prob{\left|\mathbb E_{T_n}[M_n]-\kappa\log n\right|>\delta\log n}\to 0.
  \]
  Since $\delta>0$ was arbitrary, this means that
  \[
    \frac{\mathbb E_{T_n}[M_n]}{\log n}\convprob \kappa.\qedhere
  \]
\end{proof}

We may now prove our main result.
\begin{proof}[Proof of \cref{thm:main}]
  Condition on $T_n$ and recall that
  \[
    M_n:=\sup_{v\in V(T_n)}\eta_v.
  \]
  By \cref{thm:Ding}, for any $\lambda\ge 1$,
  \[
    \left|\tau_{\mathrm{cov}}(T_n)-|E(T_n)|\left(\mathbb E_{T_n}[M_n]\right)^2\right|
    \le \lambda |E(T_n)| \sqrt{\diam(T_n)} \mathbb E_{T_n}[M_n]
  \]
  with probability at least $1-Ce^{-c\lambda}$.
  By $\diam(T_n)\le 2\cdot\height(T_n)$ and \cref{thm:height},
  \[
    \diam(T_n)
    \le
    \left(1+o_p(1)\right)\left(\frac{2}{(1+a)\xi}\right)\log n,
  \]
  and by \cref{lem:expectation},
  \[
    \mathbb E_{T_n}[M_n]=\left(1+o_p(1)\right)\frac{\gamma}{1+a}\log n.
  \]
  Since $|E(T_n)|=n-1$, the deviation term is at most
  \[
    \left(1+o_p(1)\right)\lambda n
    \left(\frac{\gamma}{1+a}\right)\left(\frac{2}{(1+a)\xi}\right)^{1/2}
    (\log n)^{3/2}
    =O_p\left(\lambda n(\log n)^{3/2}\right).
  \]
  Choose $\lambda=\lambda(n)$
  with $\lambda\to\infty$ and $\lambda=o\left((\log n)^{1/2}\right)$.
  Then $Ce^{-c\lambda}=o(1)$ and the deviation is $o_p\left(n(\log n)^2\right)$.
  Therefore,
  \[
    \tau_{\mathrm{cov}}(T_n)
    =|E(T_n)|\left(\mathbb E_{T_n}[M_n]\right)^2+o_p\left(n(\log n)^2\right)
    =\left(1+o_p(1)\right)\frac{\gamma^2}{(1+a)^2}n(\log n)^2,
  \]
  as claimed.
\end{proof}

\subsection{Expected cover time}

The following is the commute-time identity of \cite{ChandraEtAl}, specialised to trees.
\begin{lem}[Commute-time identity]\label{lem:commute-time}
  For every finite tree $T$ and every $u,v\in V(T)$,
  \[
    \mathbb E_u[\tau_v]+\mathbb E_v[\tau_u]
    =2|E(T)|d(u,v).
  \]
\end{lem}

In particular, if $T$ has $n$ vertices, then
\begin{equation}\label{eq:hit-diameter}
  t_{\mathrm{hit}}(T)\le 2(n-1)\diam(T).
\end{equation}
We will also use the following bound of Matthews.
\begin{lem}[Matthews' bound~\cite{Matthews}]\label{lem:matthews}
  For every finite connected graph $G$ on $n$ vertices,
  \[
    t_{\mathrm{cov}}(G)
    \le t_{\mathrm{hit}}(G)\sum_{k=1}^{n-1}\frac{1}{k}.
  \]
\end{lem}

\begin{proof}[Proof of \cref{cor:expectation}]
  Set $D_n:=\diam(T_n)$.
  By \eqref{eq:hit-diameter} and \cref{lem:matthews},
  \begin{equation}\label{eq:tcov-diameter}
    t_{\mathrm{cov}}(T_n)
    \le 2(n-1)D_n\sum_{k=1}^{n-1}\frac{1}{k}.
  \end{equation}
  By \eqref{eq:hit-diameter}, \eqref{eq:tcov-diameter} and \cref{thm:height},
  $t_{\mathrm{hit}}(T_n)=O_p(n\log n)$
  and $t_{\mathrm{cov}}(T_n)=O_p\left(n(\log n)^2\right)$.
  Choose $\lambda_n\to\infty$ with $\lambda_n=o((\log n)^{1/2})$;
  for example, take $\lambda_n=(\log n)^{1/4}$.
  Conditionally on $T_n$, \cref{cor:Ding} gives
  \[
    \Prob{
      \left|\tau_{\mathrm{cov}}(T_n)-t_{\mathrm{cov}}(T_n)\right|
      \ge \lambda_n
      \sqrt{t_{\mathrm{cov}}(T_n)t_{\mathrm{hit}}(T_n)}
      \,\middle|\,T_n
    }
    \le Ce^{-c\lambda_n}.
  \]
  Since $\lambda_n\sqrt{t_{\mathrm{cov}}(T_n)t_{\mathrm{hit}}(T_n)}=o_p(n(\log n)^2)$
  and $Ce^{-c\lambda_n}=o(1)$, it follows that
  \[
    \left|\tau_{\mathrm{cov}}(T_n)-t_{\mathrm{cov}}(T_n)\right|
    =o_p\left(n(\log n)^2\right).
  \]
  Together with \cref{thm:main}, this yields
  \begin{equation}\label{eq:tcov-quenched}
    \frac{t_{\mathrm{cov}}(T_n)}{n(\log n)^2}
    \convprob \frac{\gamma^2}{(1+a)^2}.
  \end{equation}

  We next justify taking expectations.
  Put $H_n:=\height(T_n)$.
  By \cref{thm:height}, $H_n/\log n$ converges in $L^1$.
  In particular,
  the family $(H_n/\log n)_{n\ge 2}$ is uniformly integrable.
  Since $D_n\le 2H_n$,
  the family $(D_n/\log n)_{n\ge 2}$ is therefore uniformly integrable.
  Moreover, \eqref{eq:tcov-diameter} gives
  \[
    \frac{t_{\mathrm{cov}}(T_n)}{n(\log n)^2}
    \le 2\frac{n-1}{n}\frac{D_n}{\log n}
    \frac{\sum_{k=1}^{n-1}1/k}{\log n}.
  \]
  The two deterministic factors on the right are bounded uniformly in $n$.
  Hence the normalised cover times form a uniformly integrable family,
  and \eqref{eq:tcov-quenched} yields
  \[
    \E{t_{\mathrm{cov}}(T_n)}
    =(1+o(1))\frac{\gamma^2}{(1+a)^2}n(\log n)^2.
    \qedhere
  \]
\end{proof}

\section*{Acknowledgements}
This note was begun at the 19th Annual Workshop on Probability and Combinatorics, held at the Bellairs Institute, Barbados, March 28--April 4, 2025.
We are very grateful to Tom Hutchcroft and James Martin for discussions at an early stage of the project.
Some further work on it took place while O.A.\ and C.G.\ were in residence at the Simons Laufer Mathematical Sciences Institute in Berkeley, California, during the Spring 2025 semester, and they gratefully acknowledge the support of the National Science Foundation under Grant No.\ DMS-1928930.
O.A.\ further thanks the Clay Institute for a Clay Senior Fellowship while at SLMath, and Magdalen College, Oxford for a visiting fellowship in Michaelmas Term 2025.

The proofs and manuscript were written by the authors.
AI tools were only used for literature research and to check the final draft for typographical and other errors.

\nocite{FuchsHwangNeininger06}
\bibliographystyle{amsplain}
\bibliography{library}

\appendix
\crefalias{section}{appendix}
\section{Biggins in the immortal immobile case}\label[appendix]{app:biggins}
In this appendix we explain why the specialised statement used in the text follows from Biggins' more general framework.

\begin{proof}[Proof of \cref{prop:biggins-right}]
  Start from the immortal branching random walk of \cref{prop:biggins-right}.
  To match Biggins' notation, take his basic triple $(Z,M,\chi)$ as follows:
  $Z$ is the offspring point process introduced in the text,
  $M(s)\equiv 0$,
  and $\chi(s)=\ind_{s\ge 0}$.
  Thus, writing $\mathcal N_t$ for Biggins' random measure (denoted $N_t$ in \cite{Biggins97}),
  $\mathcal N_t$ places one atom at the birth position of each particle that has been born by time $t$.
  Therefore, for each $x\in\R$, we have $\mathcal N_t([x,\infty))=P_t(x)$,
  and his right edge is our $R_t$.

  We now check the hypotheses of \cite[Corollary 2.2]{Biggins97}.
  First, $\chi$ is a $0$-$1$ characteristic.
  Second, by \cref{cor:alpha},
  \[
    \alpha(\theta)=a+e^{\theta^2/2}
  \]
  for every $\theta\in\R$;
  in particular $\alpha(\theta)<\infty$ for every $\theta<0$,
  which is stronger than the finiteness assumption in \cite[Corollary 2.2]{Biggins97}.
  Third, the process is well-regulated.
  By the definition just before \cite[Corollary 2.1]{Biggins97},
  this means that for each $\theta\le 0$ with $\alpha(\theta)<\infty$ one must have
  \[
    \E{\sup_{s\ge 0}
    \bigl\{
      e^{-\alpha(\theta)s}e^{-\theta M(s)}\chi(s)
    \bigr\}}<\infty.
  \]
  In our setting $M(s)\equiv 0$ and $\chi(s)=\ind_{s\ge 0}$,
  so the expression inside the supremum becomes simply $e^{-\alpha(\theta)s}$ for $s\ge 0$.
  Hence Biggins' well-regulated condition reduces to
  \[
    \sup_{s\ge 0} e^{-\alpha(\theta)s}<\infty.
  \]
  Again by \cref{cor:alpha}, $\alpha(\theta)>0$ for every $\theta$,
  so the function $s\mapsto e^{-\alpha(\theta)s}$ is decreasing on $[0,\infty)$, and therefore
  \[
    \sup_{s\ge 0} e^{-\alpha(\theta)s}=1.
  \]
  Hence the process is well-regulated.
  Thus all hypotheses of \cite[Corollary 2.2]{Biggins97} are satisfied.
  Therefore \cite[Corollary 2.2]{Biggins97} yields
  \[
    \frac{R_t}{t}\to \inf\{x:\alpha^*(x)<0\}
  \]
  almost surely.
  The quantity on the right is $\gamma_0$ by definition,
  so this is exactly the stated conclusion.
\end{proof}

\end{document}